\documentclass[12pt]{article}
\usepackage[latin1]{inputenc}
\usepackage{indentfirst}
\usepackage[dvips]{graphicx}
\usepackage{amsfonts}
\usepackage{amssymb}
\usepackage{amsmath}
\usepackage{diagrams}
\usepackage{geometry}

\begin{document}
\newtheorem{defi}{Definition}[section]
\newtheorem{prop}[defi]{Proposition}
\newtheorem{theo}[defi]{Theorem}
\newtheorem{lemma}[defi]{Lemma}
\newtheorem{ex}[defi]{Example}
\newcommand {\proof}{\emph{Proof.} }
\newcommand {\D}{\displaystyle}
\newcommand {\eps}{\varepsilon}

\title{Partial Crossed Product of a group G vs Crossed Product of S(G)}

\author{R. Exel and F. Vieira\thanks{Supported by CAPES.}}
\date{}
\maketitle

\renewcommand\abstractname{\scriptsize{ABSTRACT}}
\begin{abstract}
In this work we present a new definition to the Partial Crossed
Product by actions of inverse semigroups in a C$^*$-algebra,
without using the covariant representations as Sieben did in
\cite{sieben}. Also we present an isomorphism between the partial
crossed products by partial actions of a group $G$ and the partial
crossed product by actions of $S(G)$, an inverse semigroup
associated to $G$ introduced by Exel in \cite{exel}.
\end{abstract}

\vspace{1cm}
\section{Introduction}

In 1994 Nándor Sieben, in his master thesis, introduced the
Partial Crossed Product by an action of inverse semigroup using
covariant representations. We will present here another
definition, inspired in the definition of the partial crossed
product by a partial action. There is only one difference, we will
need to do a quotient that isn't done in the case of partial
actions. And for this, we will show that, under certain conditions
over the algebra, the Algebraic Partial Crossed Product is
associative, and this part is very similar to that Exel and
Dokuchaev done in \cite{dokuchaexel}.

Associated to a partial action $\alpha$ of a group $G$ and to the
group itself, Sieben construct a certain inverse semigroup $S_G$,
and showed that the partial crossed product of $\alpha$ is
isomorphic to the partial crossed product of a certain action of
$S_G$. Here, we present an analogous result, but using an inverse
semigroup that depends only of the group $G$. Is the semigroup
$S(G)$ introduced by Exel in \cite{exel}.

We will only work with groups $G$ with unity, that we will denote
$e$.

To start, remember that a semigroup is a set $S$ with an
associative operation.
\begin{defi}
An inverse semigroup is a semigroup $S$ such that, for any $r\in
S$, exists only one $r^*\in S$ such that $rr^*r=r$ and
$r^*rr^*=r^*$. This element $r^*$ is called inverse of $r$.
\end{defi}
\begin{ex}
For any set $X$, the set $I(X)$ of partial bijections of $X$ (that
is, bijections between subsets of $X$), is an inverse semigroup.\\
In fact, by the Wagner-Preston Theorem \cite{lawson}, any inverse
semigroup is isomorphic to an inverse subsemigroup of $I(S)$.
\end{ex}
\begin{defi}
Let $G$ be a group. We define $S(G)$ the universal semigroup
generated by $\{[g]:g\in G\}$ and by the following relations, for
$g,h\in G$:
\begin{enumerate}
    \item[(i)] $[g^{-1}][g][h]=[g^{-1}][gh]$,
    \item[(ii)] $[g][h][h^{-1}]=[gh][h^{-1}]$,
    \item[(iii)] $[g][e]=[g]$.
\end{enumerate}
\end{defi}

Note that
$[e][g]=[gg^{-1}][g]=[g][g^{-1}][g]=[g][g^{-1}g]=[g][e]=[g]$, and
then $S(G)$ is a semigroup with unit $[e]$.

The main result about this semigroup follows.
\begin{prop}\label{propu}
Let $S$ be a semigroup and $f:G\rightarrow S$ such that, for all
$g,h\in G$:
\begin{enumerate}
    \item[(i)] $f(g^{-1})f(g)f(h)=f(g^{-1})f(gh)$,
    \item[(ii)] $f(g)f(h)f(h^{-1})=f(gh)f(h^{-1})$,
    \item[(iii)] $f(g)f(e)=f(g)$.
\end{enumerate}
Then exists an unique homomorphism $\widehat{f}:S(G)\rightarrow S$
such that the diagram bellow commutes:
\begin{diagram}
G              & \rTo^f             & S\\
\dTo^{[\cdot]} & \ruTo_{\widehat{f}}&  \\
S(G)
\end{diagram}
\end{prop}
\begin{flushright}
$\square$
\end{flushright}

All details about this semigroup can be found in \cite{exel}. Some
of them we want to repeat here, because will be more important
than the others. Firstly, for $g\in G$, define
$\eps_g=[g][g^{-1}]$. Is not difficult to see that this elements
satisfy $\eps_g=\eps_g^2$ and $[h]\eps_g=\eps_{gh}[h]$, for
$g,h\in G$. Another interesting propriety is that any element
$s\in S(G)$ admits an unique decomposition
\begin{equation*}
s=\eps_{s_1}\ldots\eps_{s_n}[g],
\end{equation*}
where $n\geq0$ and $s_1,\ldots s_n,g\in G$. We will call it the
\emph{standard form} of $s$.

We can construct an anti-homomorphism $*$ in $S(G)$, such that
$[g]^*:=*([g])=[g^{-1}]$ and, for any $s\in S(G)$, $s^*$ is the
inverse (in the sense of inverse semigroups) of $s$.

With this we can conclude that that inverse is unique, and then
$S(G)$ is an inverse semigroup.

In the definition of partial crossed product by an action of
inverse semigroup that we will present, we will use a partial
order of this inverse semigroup. So, we define, for $s,t\in S$:
\begin{equation*}
s\leq t\Leftrightarrow s=tf,\hbox{ for any idempotent }f\in S.
\end{equation*}
\begin{ex}\label{exordemsg}
Take an idempotent $l=\eps_{l_1}\ldots\eps_{l_n}[j]\in S(G)$, with
$l_1,\ldots l_n,j\in G$. By the uniqueness of decomposition of
$S(G)$, is not difficult to see that $j$ must be the unit of the
group. So
$l=\eps_{l_1}\ldots\eps_{l_n}$.\\
Now, for $r=\eps_{r_1}\ldots\eps_{r_n}[h]$ and $s=\eps_{s_1}\ldots
\eps_{s_m}[g]$ in $S(G)$, if $s\leq r$ we have that $s=rf$, for
some $f\in S(G)$ idempotent. Then $f=\eps_{f_1}\ldots\eps_{f_k}$
and we have that:
\begin{equation*}
\eps_{s_1}\ldots\eps_{s_m}[g]=\eps_{r_1}\ldots\eps_{r_n}[h]\eps_{f_1}\ldots\eps_{f_k}=
\eps_{r_1}\ldots\eps_{r_n}\eps_{hf_1}\ldots\eps_{hf_k}[h].
\end{equation*}
By uniqueness of decomposition in $S(G)$, follows that $g=h$ and
$\{s_1,\ldots,s_m\}=\{r_1,\ldots,r_n,hf_1,\ldots,hf_k\}$. So, the
difference between $s$ is $r$ are some $\eps$'s.
\end{ex}
\begin{flushright}
$\square$
\end{flushright}

\section{Actions}

Let $G$ be a group with unit $e$.
\begin{defi}\label{defipaa}
A partial action $\alpha$ of $G$ in the algebra $A$ is a pair
$(\{D_g\}_{g\in G},\{\alpha_g\}_{g\in G})$ where, for each $g\in
G$, $D_g$ is a ideal of $A$ and $\alpha_g: D_{g^{-1}}\rightarrow
D_g$ is a isomorphism satisfying, for $g,h\in G$:
\begin{enumerate}
\item[(i)] $D_e=A$,
\item[(ii)] $\alpha_{g}(D_{g^{-1}}\cap
D_h)=D_g\cap D_{gh}$,
\item[(iii)]
$\alpha_g\circ\alpha_h(x)=\alpha_{gh}(x), \forall x\in
D_{h^{-1}}\cap D_{h^{-1}g^{-1}}$.
\end{enumerate}
\end{defi}

\begin{defi}\label{defipaca}
A partial action of $G$ in the C$^*$-algebra $A$ is a partial
action such that, for all $g\in G$, $D_g$ are closed ideals of $A$
and $\alpha_g$ are $*$-isomorphisms.
\end{defi}

We call $(A,G,\alpha)$ a (C$^*$-)partial dynamical system if
$\alpha$ is a partial action of the group $G$ in the
(C$^*$-)algebra $A$.

\begin{defi}\label{defiaa}
Let $S$ be an inverse semigroup with unit $e$. We say that $\beta$
is an action of $S$ in the algebra $A$ if for each $s\in S$ exists
an ideal $E_s$ of $A$ and an isomorphism $\beta_s:
E_{s^*}\rightarrow E_s$ such that, for all $r,s\in S$,
$\beta_r\circ\beta_s=\beta_{rs}$ and $E_e=A$.
\end{defi}

Is easy to verify that $(\beta_s)^{-1}=\beta_{s^*}$.

\begin{prop}\label{propuest}
Let $\beta$ an action of the inverse semigroup $S$ in the algebra
$A$. Then, for $s,t\in S$, we have that $E_{st}\subseteq E_s$.
\end{prop}
\proof By definitions follows that:
\begin{equation*}
E_{st}=\hbox{ Dom}(\beta_{t^*s^*})=\hbox{
Dom}(\beta_{t^*}\beta_{s^*})=\beta_{s^*}^{-1}(E_t\cap
E_{s^*})=\beta_s(E_t\cap E_{s^*}).
\end{equation*}
\begin{flushright}
$\square$
\end{flushright}

If the inverse semigroup is $S(G)$, we can another result:
\begin{prop}\label{propuestsu}
Let $\beta$ an action of $S(G)$ in the algebra $A$. Then, for
$g,h\in G$, we have that $E_{[g][h]}=E_{[gh]}\cap E_{[g]}$.
\end{prop}
\proof Using the last proposition, follows that:
\begin{equation*}
\begin{split}
E_{[g][h]}&=E_{[g][g^{-1}][g][h]}=E_{[g][g^{-1}][gh]}=\beta_{[g]}(E_{[g^{-1}][gh]}\cap
E_{[g^{-1}]})=\\
&=\beta_{[g]}(\beta_{[g^{-1}]}(E_{[gh]}\cap E_{[g]})\cap
E_{[g^{-1}]})=\beta_{[g]}(\beta_{[g^{-1}]}(E_{[gh]}\cap
E_{[g]}))=E_{[gh]}\cap E_{[g]}.
\end{split}
\end{equation*}
\begin{flushright}
$\square$
\end{flushright}

\begin{defi}\label{defiaca}
An action of $S$ in a $C^*$-algebra $A$ is an action of $S$ in the
algebra $A$ such that, for all $s\in S$, $E_s$ is a closed ideal
of $A$ and $\beta_s$ is a $*$-isomorphism.
\end{defi}

In \cite{exel}, Exel shows that for any group $G$ and any set $X$,
exist a bijection between the set of partial actions of $G$ on $X$
and the set of actions of $S(G)$ in a $X$.

We will extend this result, in the case that the actions act over
an algebra or a C$^*$-algebra. To start, consider the next
proposition.

\begin{prop}\label{propapaa}
Let $\beta$ an action of $S(G)$ on the algebra $A$. So
$(\{E_{[g]}\}_{g\in G},\{\beta_{[g]}\}_{g\in G})$ is a partial
action of $G$ on $A$.
\end{prop}
\proof Item $(i)$ and $(iii)$ are obvious. Item $(ii)$ follows of
$\beta_{[g]}\circ\beta_{[h]}=\beta_{[g][h]}$ and Proposition
\ref{propuestsu}.
\begin{flushright}
$\square$
\end{flushright}

The converse result follows by the universal propriety of $S(G)$.
\begin{prop}\label{proppaaa}
Let $\alpha=(\{D_g\}_{g\in G},\{\alpha_g\}_{g\in G})$ a partial
action of $G$ in the algebra $A$. So
\begin{equation*}
\begin{split}
f:G&\rightarrow I(A)\\
  g&\mapsto \alpha_g
\end{split}
\end{equation*}
satisfies $(i)$-$(iii)$ of Proposition \ref{propu}.
\end{prop}
\proof The result follows easily of items $(ii)$ and $(iii)$ of
the definition of partial action.
\begin{flushright}
$\square$
\end{flushright}

With this, we have an unique homomorphism
\begin{equation*}
\begin{split}
\beta:S(G)&\rightarrow I(A)\\
      s   &\mapsto \beta_s
\end{split}
\end{equation*}
such that $\beta_{[g]}=\alpha_g$.\\
Now, for $s_1,\ldots, s_n,h\in G$, let
$s=\eps_{s_1}\ldots\eps_{s_n}[h]=[h]\eps_{h^{-1}s_1}\ldots\eps_{h^{-1}s_n}\in
S(G)$ in the standard form. Note that
\begin{equation*}
\begin{split}
\beta_s&=\beta_{[h]}\beta_{[h^{-1}s_1]}\beta_{[(h^{-1}s_1)^{-1}]}\ldots\beta_{[(h^{-1}s_n)^{-1}]}
=\alpha_h\alpha_{h^{-1}s_1}\alpha_{(h^{-1}s_1)^{-1}}\ldots\alpha_{(h^{-1}s_n)^{-1}}=\\
&=\alpha_h|_{D_{h^{-1}s_1}\cap\ldots\cap D_{h^{-1}s_n}}.
\end{split}
\end{equation*}
So we can conclude that $E_{s^*}:=$ dom$\beta_s=D_{h^{-1}}\cap
D_{h^{-1}s_1}\cap\ldots\cap D_{h^{-1}s_n}$.\\
Since $s^*=[h^{-1}]\eps_{s_1}\ldots\eps_{s_n}$, we have $E_s:=$
dom$\beta_{s^*}=D_h\cap D_{s_1}\cap\ldots\cap D_{s_n}$.

Let us find the counter-domain of $\beta_s$. Well,
\begin{equation*}
\beta_s(E_{s^*})=\alpha_h(D_{h^{-1}}\cap
D_{h^{-1}s_1}\cap\ldots\cap
D_{h^{-1}s_n})\subseteq\alpha_h(D_{h^{-1}}\cap
D_{h^{-1}s_i})=D_h\cap D_{s_i}, \forall i.
\end{equation*}
So $\beta_s(E_{s^*})\subseteq D_h\cap D_{s_1}\cap\ldots\cap
D_{s_n}$.\\
If we change $h$ by $h^{-1}$ and $s_i$ by $h^{-1}s_i$, we conclude
that $\beta_s(E_{s^*})\supseteq D_h\cap D_{s_1}\cap\ldots\cap
D_{s_n}$, and then
\begin{equation*}
\beta_s(E_{s^*})=\alpha_h(D_{h^{-1}}\cap
D_{h^{-1}s_1}\cap\ldots\cap D_{h^{-1}s_n})=D_h\cap
D_{s_1}\cap\ldots\cap D_{s_n}=E_s.
\end{equation*}
So, for all $s\in S$, $\beta_s: E_{s^*}\rightarrow E_s$ is an
isomorphism, $E_s$ is an ideal (intersection of ideals) and like
$E_{[e]}=D_e=A$ we have that $\beta$ is an action of $S(G)$ on
$A$.

Observe that, for $r=[r_1][r_2]\ldots [r_n]\in S(G)$ (not
necessarily in the standard form),
$\beta_r=\alpha_{r_1}\alpha_{r_2}\ldots\alpha_{r_n}$,
$E_{r^*}=D_{r_n^{-1}}\cap D_{r_n^{-1}r_{n-1}^{-1}}\cap\ldots\cap
D_{r_n^{-1}r_{n-1}^{-1}\ldots r_1^{-1}}$ and $E_r=D_{r_1}\cap
D_{r_1r_2}\cap\ldots\cap D_{r_1r_2\ldots r_n}$.

So, is easy to see that holds the following theorem.
\begin{theo}\label{theobija}
Let $G$ a group and $A$ an algebra. There is a bijection between
the partial actions of $G$ on $A$ and the actions of $S(G)$ on
$A$.
\end{theo}
\begin{flushright}
$\square$
\end{flushright}

Now, if we have an action $\beta$ of $S(G)$ on the C$^*$-algebra
$A$, Proposition \ref{propapaa} implies that $(\{E_{[g]}\}_{g\in
G},\{\beta_{[g]}\}_{g\in G})$ is a partial action of $G$ on the
algebra $A$. But $E_{[g]}$ is closed and $\beta_{[g]}$ preserves
$*$. So $(\{E_{[g]}\}_{g\in G},\{\beta_{[g]}\}_{g\in G})$ is a
partial action of $G$ in the C$^*$-algebra $A$.

Conversely, if $\alpha=(\{D_g\}_{g\in G},\{\alpha_g\}_{g\in G})$
is a partial action of $G$ in the C$^*$-algebra $A$, we construct
an action $\beta$ of $S(G)$ in the algebra $A$. But intersection
of closed ideals is a closed ideal, and $\alpha_g$ preserve $*$.
So $\beta$ is an action of $S(G)$ in the C$^*$-algebra $A$. So we
can extend Theorem \ref{theobija}:
\begin{theo}\label{theobijca}
Let $G$ a group and $A$ a C$^*$-algebra. There is a bijection
between the partial actions of $G$ on $A$ and the actions of
$S(G)$ on $A$.
\end{theo}
\begin{flushright}
$\square$
\end{flushright}

\section{Algebraic Partial Crossed Product}

Consider $\alpha$ a partial action of the group $G$ (with unit
$e$) in the algebra $A$.
\begin{defi}
We define the algebraic partial crossed product of $\alpha$ as
\begin{equation*}
A\rtimes_\alpha^aG=\left\{\D\sum_{g\in
G}^{\hbox{\scriptsize{finite}}}a_g\delta_g: a_g\in D_g\right\},
\end{equation*}
where $\delta_g$ are symbols, with addition defined in the obvious
way and product being the linear extension of
\begin{equation*}
(a_g\delta_g)(a_h\delta_h)=\alpha_g(\alpha_{g^{-1}}(a_g)a_h)\delta_{gh}.
\end{equation*}
\end{defi}

Well, $\alpha_{g^{-1}}(a_g)a_h\in D_{g^{-1}}\cap D_h$, so
$\alpha_g(\alpha_{g^{-1}}(a_g)a_h)\delta_{gh}\in D_g\cap D_{gh}$,
and the multiplication is well defined.

Let us define a set that we will need later to define the
algebraic partial crossed product by an action of inverse
semigroup.
\begin{defi}\label{defiapcpasi}
Let $\beta$ an action of the inverse semigroup $S$ in the algebra
$A$. Define
\begin{equation*}
L=\left\{\D\sum_{s\in S}^{\hbox{\scriptsize{finite}}}a_s\delta_s:
a_s\in E_s\right\},
\end{equation*}
where $\delta_s$ are symbols, with addition defined in the obvious
way and product being the linear extension of
\begin{equation*}
(a_r\delta_r)(a_s\delta_s)=\beta_r(\beta_{r^{-1}}(a_r)a_s)\delta_{rs}.
\end{equation*}
\end{defi}

In \cite{dokuchaexel}, Dokuchaev and Exel proved under which
conditions $A\rtimes_\alpha G$ is associative, and they present an
example of the non-associative case (Proposition 3.6). If we take
the partial action of that Proposition, we can construct an action
of $S(G)$ (that is in bijection) such that the set $L$ above is
not associative. As we want to quotient $L$ by an ideal, we need
to know under which conditions $L$ is associative.

To this, let us talk about the algebra of multipliers.
\begin{defi}
The algebra of multipliers of a $K$-algebra $A$ is the set $M(A)$
of the ordered pairs $(L,R)$, where $L$ and $R$ are linear
transformations of $A$ such that, for $a,b\in A$:
\begin{enumerate}
\item[(i)] $L(ab)=L(a)b$,
\item[(ii)] $R(ab)=aR(b)$,
\item[(iii)]
$R(a)b=aL(b)$
\end{enumerate}
and for $(L,R)$, $(L',R')\in M(A)$, $k\in K$, the operations are
given by:
\begin{equation*}
\begin{split}
       k(L,R)&=(kL,kR),\\
(L,R)+(L',R')&=(L+L',R+R'),\\
 (L,R)(L',R')&=(L\circ L',R'\circ R).
\end{split}
\end{equation*}
We say that $L$ is a left multiplier and $R$ is a right multiplier
of $A$.
\end{defi}

We will show that that operation of multiplication $L$ is
associative when the ideals $E_s$ associated with the action
$\beta$ are $(L,R)$-associative, that is, when $L\circ R'=R'\circ
L$ for all $(L,R), (L',R')\in M(E_s)$ (for more details about the
algebra of multipliers see \cite{dokuchaexel}).
\begin{theo}
Let $\beta$ an action of the inverse semigroup $S$ in the algebra
$A$. If the ideals $E_s$ are $(L,R)$-associative, then the
operation of multiplication of the set $L$ defined in Def.
\ref{defiapcpasi} is associative.
\end{theo}
\proof Let $r$, $s$, $t\in S$ and $a_r\in E_r$, $a_s\in E_s$,
$a_t\in E_t$.\\
We want to prove that
\begin{equation*}
a_r\delta_r(a_s\delta_sa_t\delta_t)=(a_r\delta_ra_s\delta_s)a_t\delta_t,
\end{equation*}
that is:
\begin{equation*}
\beta_r(\beta_{r^*}(a_r)\beta_s(\beta_{s^*}(a_s)a_t))\delta_{rst}=
\beta_{rs}(\beta_{s^*r^*}(\beta_r(\beta_{r^*}(a_r)a_s))a_t)\delta_{rst}.
\end{equation*}
Analyzing the right side:
\begin{equation*}
\begin{split}
\beta_{rs}(\beta_{s^*r^*}(\beta_r(\beta_{r^*}(a_r)a_s))a_t)\delta_{rst}&=
\beta_{rs}(\beta_{s^*}(\beta_{r^*}(\beta_r(\beta_{r^*}(a_r)a_s)))a_t)\delta_{rst}=\\
&=\beta_{rs}(\beta_{s^*}(\beta_{r^*}(a_r)a_s)a_t)\delta_{rst}=\\
&=\beta_r(\beta_s(\beta_{s^*}(\beta_{r^*}(a_r)a_s)a_t))\delta_{rst}.
\end{split}
\end{equation*}
So we need to prove that:
\begin{equation*}
\beta_r(\beta_{r^*}(a_r)\beta_s(\beta_{s^*}(a_s)a_t))=
\beta_r(\beta_s(\beta_{s^*}(\beta_{r^*}(a_r)a_s)a_t)).
\end{equation*}
Applying $\beta_{r^*}$ in both sides of equality by the left:
\begin{equation*}
\beta_{r^*}(\beta_r(\beta_{r^*}(a_r)\beta_s(\beta_{s^*}(a_s)a_t)))=
\beta_{r^*}(\beta_r(\beta_s(\beta_{s^*}(\beta_{r^*}(a_r)a_s)a_t)))
\end{equation*}
and it is equivalent to:
\begin{equation*}
\beta_{r^*}(a_r)\beta_s(\beta_{s^*}(a_s)a_t)=
\beta_s(\beta_{s^*}(\beta_{r^*}(a_r)a_s)a_t).
\end{equation*}
Because $\beta_{r^*}:E_r\rightarrow E_{r^*}$ is an isomorphism,
the above condition is equivalent to:
\begin{equation*}
a\beta_s(\beta_{s^*}(a_s)a_t)=
\beta_s(\beta_{s^*}(aa_s)a_t),\hbox{ }\forall a\in E_{r^*},\hbox{
}a_s\in E_s,\hbox{ }a_t\in E_t.
\end{equation*}
Denoting $R_{a_t}:E_{s^*}\rightarrow E_{s^*}$ the right multiplier
by $a_t$ in $E_{s^*}$ and $L_a:E_s\rightarrow E_s$ the left
multiplier by $a$ in $E_s$, the last equation is equivalent to:
\begin{equation*}
L_a\circ\beta_s\circ R_{a_t}\circ\beta_{s^*}(a_s)=\beta_s\circ
R_{a_t}\circ\beta_{s^*}\circ L_a(a_s),\hbox{ }\forall a\in
E_{r^*},\hbox{ }a_s\in E_s,\hbox{ }a_t\in E_t.
\end{equation*}
Now, $\beta_s\circ R_{a_t}\circ\beta_{s^*}$ is a right multiplier
of $E_s$, and because $E_s$ is $(L,R)$-associative, the last
equation holds.\\
So the multiplication of $L$ is associative.
\begin{flushright}
$\square$
\end{flushright}

So, let us suppose that the ideals $E_s$ related with the action
$\beta$ of $S$ in $A$ are $(L,R)$-associative.
\begin{defi}
Let $\beta$ an action of the inverse semigroup $S$ in the algebra
$A$. Consider $N=\langle a\delta_r-a\delta_t:a\in E_r, r\leq
t\rangle$, that is, the ideal generated by $a\delta_r-a\delta_t$.
We define the algebraic partial crossed product of $\beta$ as
\begin{equation*}
A\rtimes_\beta^a S=\frac{L}{N}.
\end{equation*}
\end{defi}

Note that $r\leq t$ implies $r=ti$, for $i$ idempotent. So, by
Proposition \ref{propuest}, $E_r\subseteq E_t$.

Because the definition, we denote the elements of
$A\rtimes_\beta^a S$ like $\overline{a_s\delta_s}$, where
$a_s\delta_s\in L$.
\begin{lemma}\label{lemma36}
Let $\beta$ an action of $S$ in $A$. For $r_1,\ldots,r_n,g,h\in
G$, hold in $A\rtimes_\beta^a S$:
\begin{enumerate}
    \item[(1)]
    $\overline{a\delta_{[g][h]}}=\overline{a\delta_{[gh]}}$, for
    $a\in E_{[g][h]}$,
    \item[(2)]
    $\overline{a\delta_{\eps_{r_1}\ldots\eps_{r_n}[g]}}=\overline{a\delta_{[g]}}$,
    for $a\in E_{\eps_{r_1}\ldots\eps_{r_n}[g]}$.
\end{enumerate}
\end{lemma}
\proof $(1)$: Well, $[g][h]=[g][h][h^{-1}][h]=[gh][h^{-1}][h]$. As
$[h^{-1}][h]$ is idempotent, $[g][h]\leq [gh]$ and so
$a\delta_{[g][h]}-a\delta_{[gh]}\in N$.\\
$(2)$: Note that $\eps_{r_1}\ldots\eps_{r_n}[g]=
[g]\eps_{g^{-1}r_1}\ldots\eps_{g^{-1}r_n}$ and the results follows
because $\eps_{g^{-1}r_1}\ldots\eps_{g^{-1}r_n}$ is idempotent.
\begin{flushright}
$\square$
\end{flushright}

So, now we can enunciate the main result of this section:
\begin{theo}\label{theoisoa}
Let $\alpha$ a partial action of the group $G$ in the algebra $A$.
Consider the action $\beta$ related with $\alpha$ by the Theo.
\ref{theobija}. Then $A\rtimes_\alpha^aG\cong
A\rtimes_\beta^aS(G)$.
\end{theo}
\proof Define
\begin{equation*}
\begin{split}
\varphi:A\rtimes_\alpha^a G&\rightarrow A\rtimes_\beta^a
S(G)\\
a\delta_g&\mapsto\overline{a\delta_{[g]}},\hbox{ linearly
extended.}
\end{split}
\end{equation*}
Let us prove that $\varphi$ is an isomorphism. It is well defined
and using Lemma above we can prove that it is a homomorphism,
because:
\begin{equation*}
\begin{split}
\varphi(a\delta_g)\varphi(b\delta_h)&=(\overline{a\delta_{[g]}})
(\overline{b\delta_{[h]}})=\overline{\beta_{[g]}(\beta_{[g^{-1}]}(a)b)\delta_{[g][h]}}=
\overline{\beta_{[g]}(\beta_{[g^{-1}]}(a)b)\delta_{[gh]}},\\
\varphi((a\delta_g)(b\delta_h))&=\varphi(\alpha_g(\alpha_{g^{-1}}(a)b)\delta_{gh})=
\overline{\alpha_g(\alpha_{g^{-1}}(a)b)\delta_{[gh]}}=
\overline{\beta_{[g]}(\beta_{[g^{-1}]}(a)b)\delta_{[gh]}}.
\end{split}
\end{equation*}
To show that $\varphi$ is bijective, let present an inverse for
it. For $s_1,\ldots s_n,g\in G$ consider $s=\eps_{s_1}\ldots
\eps_{s_n}[g]\in S(G)$ and the function $\gamma$ such that
$\gamma(s)=g$. Is very easy to show that $\gamma$ is a
homomorphism between the inverses semigroups $S(G)$ and $G$. Note
that $\gamma(s^*)=g^{-1}=\gamma(s)^{-1}$. So define:
\begin{equation*}
\begin{split}
\psi:L&\rightarrow A\rtimes_\alpha^a G\\
a\delta_s&\mapsto a\delta_{\gamma(s)},\hbox{ linearly extended.}
\end{split}
\end{equation*}
Note that $\psi$ is a homomorphism and using Example
\ref{exordemsg}
we see that $\psi(N)=0$.\\
So we can extend $\psi$ to the homomorphism
\begin{equation*}
\begin{split}
\widetilde{\psi}:A&\rtimes_\beta^a S(G)\rightarrow
A\rtimes_\alpha^a
G\\
\overline{a\delta_s}&\mapsto a\delta_{\gamma(s)}\hbox{, linearly
extended.}
\end{split}
\end{equation*}
As is obvious that $\widetilde{\psi}$ and $\varphi$ are inverses
one each other, the theorem holds.
\begin{flushright}
$\square$
\end{flushright}

\section{Partial Crossed Product}

Let $A$ a C$^*$-algebra with unit. For $g\in G$ and $a_g\in D_g$,
define in $A\rtimes_\alpha^aG$ the following operation $^*$:
\begin{equation*}
\begin{split}
(a_g\delta_g)^*     &=\alpha_{g^{-1}}(a_g^*)\delta_{g^{-1}},\hbox{ linearly extended:}\\
\left(\D\sum_{g\in G}^{\hbox{\scriptsize{finite}}}
a_g\delta_g\right)^*&=\D\sum_{g\in
G}^{\hbox{\scriptsize{finite}}}(a_g\delta_g)^*.
\end{split}
\end{equation*}

Is easy to show that $A\rtimes_\alpha^aG$ with $^*$ is a
$*$-algebra. Considering the following norm in
$A\rtimes_\alpha^aG$, we have that it is a normed $*$-algebra:
\begin{equation*}
\left\|\D\sum_{g\in
G}^{\hbox{\scriptsize{finite}}}a_g\delta_g\right\|_1=\D\sum_{g\in
G}^{\hbox{\scriptsize{finite}}}\|a_g\|,
\end{equation*}
where the norm in the right side is the norm in $A$.

Given a Banach $*$-algebra $B$, its enveloping C$^*$-algebra is
the completion of $B/$ker$\rho_s$ with respect to $\rho_s(x)=
\hbox{sup}\{\|\rho(x)\|:\rho\hbox{ is a representation of }B\}$.

To define the partial crossed product by a partial action $\alpha$
of $G$ in $A$, we want to take the enveloping C$^*$-algebra of
$A\rtimes_\alpha^aG$. But this set isn't a Banach $*$-algebra. So
we need to show that its representations are contractive. To this,
let $\pi$ a representation of $A\rtimes_\alpha^aG$. For $a_g\in
D_g$:
\begin{equation*}
\|\pi(a_g\delta_g)\|^2=\|\pi(a_g\delta_g)^*\pi(a_g\delta_g)\|=
\|\pi((a_g\delta_g)^*(a_g\delta_g))\|=\|\pi(\alpha_{g^{-1}}(a_g^*a_g)\delta_e)\|.
\end{equation*}
Now, note that $\alpha_{g^{-1}}(a_g^*a_g)\delta_e\in A\delta_e$,
and $A\delta_e$ is a C$^*$-algebra (isomorphic to $A$). Then:
\begin{equation*}
\|\pi(a_g\delta_g)\|^2=\|\pi(\alpha_{g^{-1}}(a_g^*a_g)\delta_e)\|\leq
\|\alpha_{g^{-1}}(a_g^*a_g)\delta_e\|_1=\|\alpha_{g^{-1}}(a_g^*a_g)\|=
\|a_g^*a_g\|=\|a_g\|^2.
\end{equation*}
So, we have that:
\begin{equation*}
\left\|\pi\left(\D\sum_{g\in
G}^{\hbox{\scriptsize{finite}}}a_g\delta_g\right)\right\|\leq\D\sum_{g\in
G}^{\hbox{\scriptsize{finite}}}\|\pi(a_g\delta_g)\|\leq\D\sum_{g\in
G}^{\hbox{\scriptsize{finite}}}\|a_g\|=\left\|\D\sum_{g\in
G}^{\hbox{\scriptsize{finite}}}a_g\delta_g\right\|_1.
\end{equation*}

Then we can take the enveloping C$^*$-algebra of
$A\rtimes_\alpha^aG$ to define the partial crossed product.
\begin{defi}
The partial crossed product of the group of the partial action
$\alpha$ of $G$ in $A$, denoted $A\rtimes_\alpha G$, is the
enveloping C$^*$-algebra of the $*$-algebra $A\rtimes_\alpha^aG$.
\end{defi}

Let us denote the elements of $A\rtimes_\alpha G$ as classes of
the elements of $A\rtimes_\alpha^aG$, $\overline{a_g\delta_g}$.

To define the partial crossed product by an action $\beta$ of an
inverse semigroup $S$ in $A$, we want to do the same, that is,
take the enveloping C$^*$-algebra of $A\rtimes_\alpha^aS$.

So, for $r\in S$ e $a_r\in E_r$, define in $A\rtimes_\beta^aS$:
\begin{equation*}
\left(\overline{a_r\delta_r}\right)^*=\overline{\beta_{r^*}(a_r^*)\delta_{r^*}},\hbox{
linearly extended.}
\end{equation*}

Easily we can see that $A\rtimes_\beta^aS$ is a $*$-algebra. Also
define a norm
\begin{equation*}
\left\|\D\sum_{s\in
S}^{\hbox{\scriptsize{finite}}}a_s\delta_s\right\|_2=\D\sum_{s\in
S}^{\hbox{\scriptsize{finite}}}\|a_s\|,
\end{equation*}
and is easy to check that $A\rtimes_\beta^aS$ is a normed
$*$-algebra.
\begin{prop}\label{propconts}
Every representation $\rho:A\rtimes_\beta^aS\rightarrow B(H)$ is
contractive.
\end{prop}
\proof A representation $\rho$ of $A\rtimes_\beta^aS$ is one of
$L$ such that $\rho|_N\equiv0$.\\
Let $s\in S$ and $a_s\in E_s$. As $s^*s\leq e$ ($e$ the unit of
$S$), follows that
\begin{equation*}
\begin{split}
\|\rho(a_s\delta_s)\|^2&=\|\rho(a_s\delta_s)^*\rho(a_s\delta_s)\|=
\|\rho((a_s\delta_s)^*(a_s\delta_s))\|=\|\rho(\beta_{s^*}(a_s^*a_s)\delta_{s^*s})\|=\\
&=\|\rho(\beta_{s^*}(a_s^*a_s)\delta_{e})\|\leq\|a_s^*a_s\|=\|a_s\|^2.
\end{split}
\end{equation*}
For any element $\D\sum_{s\in
S}^{\hbox{\scriptsize{finite}}}a_s\delta_s\in L$:
\begin{equation*}
\begin{split}
\left\|\rho\left(\D\sum_{s\in
S}^{\hbox{\scriptsize{finite}}}a_s\delta_s\right)\right\|&\leq
\D\sum_{s\in
S}^{\hbox{\scriptsize{finite}}}\|\rho(a_s\delta_s)\|\leq
\D\sum_{s\in
S}^{\hbox{\scriptsize{finite}}}\|a_s\|=\left\|\D\sum_{s\in
S}^{\hbox{\scriptsize{finite}}}a_s\delta_s\right\|_2.
\end{split}
\end{equation*}
So $\rho$ is contractive.
\begin{flushright}
$\square$
\end{flushright}

\begin{defi}: \label{defipcps}
The partial crossed product by the action $\beta$ of the inverse
semigroup $S$ on the C$^*$-algebra $A$, denoted $A\rtimes_\beta
S$, is the enveloping C$^*$-algebra of the $*$-algebra
$A\rtimes_\beta^aS$.
\end{defi}

Note that to construct $A\rtimes_\beta S$, we do two quotients.
So, we will denote its elements like
$\overline{\overline{a_r\delta_r}}$, where $a_r\delta_r\in L$.

Let $\alpha$ a partial action of the group $G$ in the
C$^*$-algebra $A$. Consider the action $\beta$ of $S(G)$ in $A$
related by the Theorem \ref{theobijca}.
\begin{theo}
The C$^*$-algebras $A\rtimes_\alpha G$ and $A\rtimes_\beta S$ are
isomorphic.
\end{theo}
\proof For $g\in G$ e $a_g\in D_g$, define
\begin{equation*}
\begin{split}
\phi: A\rtimes_\alpha^aG&\rightarrow A\rtimes_\beta S(G)\\
             a_g\delta_g&\rightarrow
             \overline{\overline{a_g\delta_{[g]}}}\hbox{, linearly extended.}
\end{split}
\end{equation*}
Obviously $\phi$ is well defined and, using Lemma \ref{lemma36},
we see that is a homomorphism. Also is easy to check that $\phi$
preserves $^*$.\\
Then, by the enveloping propriety of $A\rtimes_\alpha G$, follows
that exists unique $*$-homomorphism $\varphi: A\rtimes_\alpha
G\rightarrow A\rtimes_\beta S(G)$ such that the diagram bellow
commutes, that is, $\varphi(\overline{a_g\delta_g})=
\overline{\overline{a_g\delta_{[g]}}}$.
\begin{diagram}
A\rtimes_\alpha^aG  &   \rTo^\phi       &   A\rtimes_\beta S(G)\\
\dTo^{\iota}        &   \ruTo^{\varphi} &    \\
A\rtimes_\alpha G
\end{diagram}
Let \begin{center} $K=\left\{\displaystyle\sum_{s\in
S(G)}^{\hbox{\scriptsize{finito}}}a_s\delta_s:a_s\in
E_s\right\}$.\end{center} Using the homomorphism
$\gamma:S(G)\rightarrow G$ that we define in Theorem
\ref{theoisoa}, consider
\begin{equation*}
\begin{split}
\omega:    K&\rightarrow A\rtimes_\alpha G\\
 a_s\delta_s&\mapsto
             \overline{a_s\delta_{\gamma(s)}}\hbox{, linearly extended.}
\end{split}
\end{equation*}
Now take $M$ the ideal of $K$ generated by $a\delta_r-a\delta_t$,
where $a\in E_r$ and $r\leq t$. Is easy to see that $\omega$ is a
homomorphism and, using Example \ref{exordemsg}, also we can
see that $\omega(M)=0$.\\
So we can define
\begin{equation*}
\begin{split}
\widetilde{\phi}: A\rtimes_\beta^a S(G)&\rightarrow A\rtimes_\alpha G\\
                 \overline{a_s\delta_s}&\mapsto
             \overline{a_s\delta_{\gamma(s)}}\hbox{, linearly extended.}
\end{split}
\end{equation*}
By Lemma \ref{lemma36}, $\widetilde{\phi}$ is a homomorphism and
easily is checked that preserves $^*$. Then, by the enveloping
propriety of $A\rtimes_\beta S(G)$ exists unique $*$-homomorphism
$\widetilde{\varphi}: A\rtimes_\beta S(G)\rightarrow
A\rtimes_\alpha G$ such that the diagram bellow commutes,
\begin{diagram}
A\rtimes_\beta^a S(G) &   \rTo^{\widetilde{\phi}}     &   A\rtimes_\alpha G\\
\dTo^{\iota}          &   \ruTo^{\widetilde{\varphi}} &    \\
A\rtimes_\beta S(G)
\end{diagram}
that is, $\widetilde{\varphi}(\overline{\overline{a_s\delta_s}})=
\overline{a_s\delta_{\gamma(s)}}$.\\
Obviously $\varphi\circ\widetilde{\varphi}=Id_{A\rtimes_\beta
S(G)}$ and $\widetilde{\varphi}\circ\varphi=Id_{A\rtimes_\alpha
G}$, and the theorem holds.
\begin{flushright}
$\square$
\end{flushright}

\section{Covariant Representations}

The partial crossed product by an action of inverse semigroup was
introduced by Nándor Sieben in \cite{sieben} in 1994. In this
definition, he used covariant representations of an action. Here
we will present them and show that our definition is equivalent to
that he present. To start, is very straightforward to show the
next Proposition.
\begin{prop}
Let $H, K$ Hilbert spaces. The following conditions on an operator
$U\in B(H,K)$ are equivalent:
\begin{enumerate}
    \item[(1)] $U=UU^*U$,
    \item[(2)] $P=U^*U$ is a projection,
    \item[(3)] $U|_{\hbox{ker}^\perp U}$ is an isometry.
\end{enumerate}
An operator which satisfies the equivalent conditions above is
called a partial isometry. Let us denote $PIso(H,K)$ the set of
partial isometries between $H$ and $K$.
\end{prop}
\begin{flushright}
$\square$
\end{flushright}

Let $U\in PIso(H,K)$. If we define $M=ker^\perp(U)$ and $N=Im(U)$
the image of $U$, we conclude that $U\equiv 0$ in $M^\perp$ and
that besides being surjective, $U$ is an unitary isometry between
$M$ and $N$. In fact:
\begin{equation*}
UU^*(N)=UU^*U(M)=U(M)=N.
\end{equation*}
So $UU^*=Id_N$.

We call $M$ the \emph{initial space} and $N$ the \emph{final
space} of $U$.

For more details about the partial isometries, see section 4.2 of
\cite{sunder}.
%\begin{defi}
%Let $G$ be a group with unit $e$ and $H$ a Hilbert space. A
%partial representation of $G$ in $H$ is an application $\rho:
%G\rightarrow B(H)$ such that, for $g,h\in G$:
%\begin{enumerate}
%\item[(i)] $\rho(g)\rho(h)\rho(h^{-1})=\rho(gh)\rho(h^{-1})$,
%\item[(ii)] $\rho(g^{-1})=\rho(g)^*$,
%\item[(iii)] $\rho(e)=Id_H$.
%\end{enumerate}
%\end{defi}

So, we can define a covariant representation of an action of
inverse semigroup.
\begin{defi}
Let $\beta$ an action of the inverse semigroup $S$ in the
C$^*$-algebra $A$. A covariant representation of $\beta$ is a
triple $(\pi,\nu,H)$ where $\pi: A\rightarrow B(H)$ is a
representation of $A$ in the Hilbert space $H$ and $\nu:
S\rightarrow PIso(H)$ preserves product such that, for $s\in S$:
\begin{enumerate}
    \item[(i)] $\nu_s\pi(a)\nu_{s^*}=\pi(\beta_s(a))$ for all $a\in E_{s^*}$ (covariance condition),
    \item[(ii)] $\nu_s$ has initial space $\overline{\hbox{span}\{\pi(E_{s^*})H\}}$
    and final space $\overline{\hbox{span}\{\pi(E_s)H\}}$.
\end{enumerate}
\end{defi}

The set of covariant representations of $(A,S,\beta)$ is denoted
CovRep$(A,S,\beta)$.

Now, let us define the partial crossed product like Sieben do in
his work \cite{sieben}. All demonstrations can be founded in that
article.

Let $\beta$ an action of the unital inverse semigroup $S$ in the
C$^*$-algebra $A$. Define
\begin{center}
$\widetilde{L}=\left\{x\in l^1(S,A):x(s)\in E_s\right\}$,
\end{center}
with norm, scalar multiplication and addition inherited of
$l^1(S,A)$.\\
To $x,y\in\widetilde{L}$, the product $x*y$ is defined:
\begin{equation*}
(x*y)(s)=\displaystyle\sum_{rt=s}\beta_r(\beta_{r^*}(x(r))y(t)).
\end{equation*}
Beyond this, define $x^*$ the element of $l^1(S,A)$ such that:
\begin{equation*}
x^*(s)=\beta_s(x(s^*)^*).
\end{equation*}

The operations are well defined and $\widetilde{L}$ is a Banach
$*$-algebra.

\begin{defi}
If $(\pi,\nu,H)\in$ CovRep$(A,S,\beta)$, e define $\pi\times\nu:
\widetilde{L}\rightarrow B(H)$ as
\begin{center}
$(\pi\times\nu)(x)=\displaystyle\sum_{s\in S}\pi(x(s))\nu_s$.
\end{center}
\end{defi}

We have that $\pi\times\nu$ is a $*$-homomorphism.

Nándor Sieben defines the partial crossed product by an action of
inverse semigroup as follows.
\begin{defi}
Let $\beta$ an action of the unital inverse semigroup $S$ in the
C$^*$-algebra $A$. Define a seminorm $\|.\|_c$ in $\widetilde{L}$
as
\begin{center}
$\|x\|_c=\sup\{\|(\pi\times\nu)(x)\|:
(\pi,\nu,H)\in\hbox{CovRep}(A,S,\beta)\}$.
\end{center}
Consider $I=\{x\in\widetilde{L}:\|x\|_c=0\}$. The partial crossed
product of $\beta$ is the C$^*$-algebra gotten by the completion
of the quotient $\widetilde{L}/I$ with respect to $\|.\|_c$.
\end{defi}

We want to prove that the above definition is equivalent to
Definition \ref{defipcps}. The first step in this way is to show
that in the above definition, we can take the set $L$ instead of
$\widetilde{L}$.
\begin{lemma}
Let $E\subseteq F$ linear spaces and $\|.\|_1$ e $\|.\|_2$ two
norms over $F$ such that exist $k$(constant) such that, for all
$x\in F$, $\|x\|_2\leq k\|x\|_1$. Suppose that $E$ is dense in $F$
with respect to $\|.\|_1$. Then, the completions of $E$ and $F$
with respect to $\|.\|_2$ are isomorphic.
\end{lemma}
\proof For $i=1,2$ denote $\overline{E}^i$ the completion of $E$
with respect to $\|.\|_i$, and the same for $F$. The function
$T:(E,\|.\|_2)\rightarrow\overline{F}^2$, that includes an element
of $E$ in $F$ and then put them in $\overline{F}^2$, is a linear
isometry, and then is uniformly continuous. So we can extend it to
the isometry $\widetilde{T}: \overline{E}^2
\rightarrow\overline{F}^2$. Is easy to see that $F\subseteq
Im(\widetilde{T})$.\\
Now, as $\widetilde{T}$ is an isometry and $\overline{E}^2$ is
complete, $Im(\widetilde{T})$ is complete with respect to
$\|.\|_2$. Then, $Im(\widetilde{T})=\overline{F}^2$ and
$\widetilde{T}$ is an isomorphism between $\overline{E}^2$ and
$\overline{F}^2$.
\begin{flushright}
$\square$
\end{flushright}

So, consider $L\subseteq\widetilde{L}$. In $\widetilde{L}$ we have
defined a norm
\begin{equation*}
\|x\|_a=\D\sum_{s\in S}\|x(s)\|,
\end{equation*}
and a seminorm
\begin{equation*}
\|x\|_c=\sup\{\|(\pi\times\nu)(x)\|:
(\pi,\nu,H)\in\hbox{CovRep}(A,S,\beta)\}.
\end{equation*}

Note that we can take the completion of $\widetilde{L}$ with
respect to the seminorm $\|.\|_c$, and this is equal to the
completion of $\widetilde{L}/(\hbox{ker}\|.\|_c)$ with respect to
the norm $\|.\|_c$.

Also note that, for $x\in\widetilde{L}$ and
$(\pi,\nu,H)\in\hbox{CovRep}(A,S,\beta)$:
\begin{equation*}
\begin{split}
\|(\pi\times\nu)(x)\|&=\left\|\D\sum_{s\in
S}\pi(x(s))\nu_s\right\|\leq \D\sum_{s\in S}\|\pi(x(s))\nu_s\|
\leq \D\sum_{s\in S}\|\pi(x(s))\|\|\nu_s\|\leq\\
&\leq\D\sum_{s\in S}\|\pi(x(s))\|\leq\D\sum_{s\in
S}\|x(s)\|=\|x\|_a,
\end{split}
\end{equation*}
and then $\|x\|_c=\sup\|(\pi\times\nu)(x)\|\leq\|x\|_a$. As
$\widetilde{L}$ is the completion of $L$ with respect to $\|.\|_a$
follows, by the previous Lemma, that in the definition that Sieben
present, we can consider the set $L$ instead of $\widetilde{L}$.

Observe that by the Cohen-Hewitt Factorization Theorem (Theorem
32.22, \cite{hewitt}), is easy to show that
$\overline{span\{\pi(I)H\}}=\pi(I)H$, for any $I$ closed ideal of
the C$^*$-algebra $A$.

So, if we show that
\begin{equation*}
\|x\|_c=\rho_s(x),
\end{equation*}
for all $x\in L$, we have proved that the two definitions of
partial crossed product by an action of inverse semigroup are the
same. To this, consider the next theorem.
\begin{theo}
Let $\rho$ be a representation of $L$ in the Hilbert space $H$.
Then $\rho|_N\equiv 0\Leftrightarrow\rho=\pi\times\nu$ for some
$(\pi,\nu,H)\in$ CovRep$(A,S,\beta)$.
\end{theo}
\proof $(\Leftarrow)$ Let $\rho=\pi\times\nu$ and take
$a\delta_r-a\delta_t$ a generator of $N$, that is, $r=ti$, for $i$
idempotent. For $a\in E_r=E_{it}\subseteq E_i$ and $h\in H$:
\begin{equation*}
\begin{split}
\rho(a\delta_r-a\delta_t)(h)&=(\pi\times\nu)(a\delta_r-a\delta_t)(h)=\pi(a)\nu_r(h)-\pi(a)\nu_t(h)=\\
&=\pi(a)\nu_{it}(h)-\pi(a)\nu_t(h)=\pi(a)\nu_i\nu_t(h)-\pi(a)\nu_t(h).
\end{split}
\end{equation*}
Denote $K_i=\pi(E_i)H$. Then $H=K_i\oplus(K_i)^\perp$, and let us
do the demonstration in cases:

\underline{$\nu_t(h)\in K_i$:} As $i$ is idempotent, follows that
it is the identity in $K_i$. So $\nu_i(\nu_t(h))=\nu_t(h)$ and
then:
\begin{equation*}
\begin{split}
\pi(a)\nu_i\nu_t(h)-\pi(a)\nu_t(h)=\pi(a)\nu_t(h)-\pi(a)\nu_t(h)=0.
\end{split}
\end{equation*}

\underline{$\nu_t(h)\in (K_i)^\perp$:} Note that $\pi(E_i)\equiv
0$ in $(K_i)^\perp$ and then:
\begin{equation*}
\begin{split}
\pi(a)\nu_i\nu_t(h)-\pi(a)\nu_t(h)=0.
\end{split}
\end{equation*}
So $\rho(N)=0$.\\
$(\Rightarrow)$ Suppose that $\rho|_N\equiv0$. So, for $r\leq t$
and $a\in E_r$, $\rho(a\delta_r)=\rho(a\delta_t)$. Define\\
\begin{minipage}{.4\linewidth}
\begin{equation*}
\begin{split}
\pi:A&\rightarrow B(H)\\
    a&\mapsto\rho(a\delta_e),
\end{split}
\end{equation*}
\end{minipage}
\begin{minipage}{.5\linewidth}
\begin{equation*}
\begin{split}
\nu:S&\rightarrow B(H)\\
    s&\mapsto\D\lim_\lambda\rho(u_\lambda\delta_s),
             \{u_\lambda\}\hbox{ approx. identity of }E_s,
\end{split}
\end{equation*}
\end{minipage}
\\
\\
and $\D\lim_\lambda$ denotes the strong operator limit.\\
Let us prove that $(\pi,\nu,H)\in$ CovRep$(A,S,\beta)$. Is obvious
that $\pi$ is a representation.\\
To show that $\nu$ is well defined, let $s\in S$ and consider
$\{u_\lambda\}$ an approximately identity for $E_s$. As
$\beta_{s^*}:E_s\rightarrow E_{s^*}$ is an isomorphism, we know
that $\{\beta_{s^*}(u_\lambda)\}$ is an approximately identity of
$E_{s^*}$. As $H=\pi(E_{s^*})H\oplus(\pi(E_{s^*})H)^\perp$, we
will split the demonstration. If $h\in\pi(E_{s^*})H$, then
$h=\pi(a)k=\rho(a\delta_e)k$, for $a\in E_{s^*},k\in H$. So:
\begin{equation*}
\begin{split}
\nu_s(h)&=\D\lim_\lambda\rho(u_\lambda\delta_s)(h)=\D\lim_\lambda\rho(u_\lambda\delta_s)\rho(a\delta_e)(k)=
\D\lim_\lambda\rho(\beta_s(\beta_{s^*}(u_\lambda)a)\delta_s)(k)=\\
&=\rho(\beta_s(a)\delta_s)(k).
\end{split}
\end{equation*}
If $h\in(\pi(E_{s^*})H)^\perp$, we have that $\langle
h,\rho(E_{s^*}\delta_e)H\rangle=\langle h,\pi(E_{s^*})H\rangle=0$.
So:
\begin{equation*}
\begin{split}
\langle\rho(\beta_{s^*}(\sqrt{u_\lambda})\delta_e)(h),H\rangle=\langle
h,\rho(\beta_{s^*}(\sqrt{u_\lambda})\delta_e)H\rangle=0,
\end{split}
\end{equation*}
that implies $\rho(\beta_{s^*}(\sqrt{u_\lambda})\delta_e)(h)=0$.
Then:
\begin{equation*}
\begin{split}
\D\lim_\lambda\rho(u_\lambda\delta_s)(h)&=
\D\lim_\lambda\rho[\beta_s(\beta_{s^*}(\sqrt{u_\lambda})\beta_{s^*}(\sqrt{u_\lambda}))\delta_s](h)=\\
&=\D\lim_\lambda\rho[(\sqrt{u_\lambda}\delta_s)(\beta_{s^*}(\sqrt{u_\lambda})\delta_e)](h)=\\
&=\D\lim_\lambda\rho(\sqrt{u_\lambda}\delta_s)\rho(\beta_{s^*}(\sqrt{u_\lambda})\delta_e)(h)=0.
\end{split}
\end{equation*}
So $\nu_s$ is independent of the approximately identity taken. As
$\rho$ is contractive (Proposition \ref{propconts}):
\begin{equation*}
\begin{split}
\|\nu_s\|=\|\D\lim_\lambda\rho(u_\lambda\delta_s)\|=\D\lim_\lambda\|\rho(u_\lambda\delta_s)\|\leq
\D\lim_\lambda\|u_\lambda\delta_s\|\leq\D\lim_\lambda\|u_\lambda\|\leq
1,
\end{split}
\end{equation*}
and then $\nu_s\in B(H)$. So it is well defined.\\
To show that $\nu_s$ is a partial isometry (with initial space
$\pi(E_{s^*})H$ and final $\pi(E_s)H$), first let us show that
$\nu_s^*=\nu_{s^*}$. Let $\{u_\lambda\}$ an approximately identity
of $E_{s^*}$. Then, for $k_1,k_2\in H$:
\begin{equation*}
\begin{split}
\langle k_1,\nu_{s^*}(k_2)\rangle&=\langle
k_1,\D\lim_\lambda\rho(u_\lambda\delta_{s^*})(k_2)\rangle=\D\lim_\lambda\langle
k_1,\rho(u_\lambda\delta_{s^*})(k_2)\rangle=\D\lim_\lambda\langle\rho(u_\lambda\delta_{s^*})^*(k_1),k_2\rangle=\\
&=\langle\D\lim_\lambda\rho(\beta_s(u_\lambda)\delta_s)(k_1),k_2\rangle=\langle\nu_s(k_1),k_2\rangle.
\end{split}
\end{equation*}
So $\nu_s^*=\nu_{s^*}$.\\
Let us show that $\nu_s^*\nu_s$ is a projection over
$\pi(E_{s^*})H$, because we see yet that $\nu_s\equiv 0$ in
$(\pi(E_{s^*})H)^\perp$. Then, let $h=\pi(a)k\in\pi(E_{s^*})H$ and
$\{u_\gamma\}$ an approximately identity of $E_{s^*}$:
\begin{equation*}
\begin{split}
\nu_s^*\nu_s(h)&=\nu_{s^*}(\rho(\beta_s(a)\delta_s)(k))=
\D\lim_\gamma\rho(u_\gamma\delta_{s^*})\rho(\beta_s(a)\delta_s)(k)=\\
&=\D\lim_\gamma\rho(\beta_{s^*}(\beta_s(u_\gamma)\beta_s(a))\delta_{s^*s})(k)=
\rho(a\delta_{s^*s})(k)=\\
&=\rho(a\delta_e)(k)=h,
\end{split}
\end{equation*}
because $\rho|_N=0$. Then $\nu_s$ is a partial isometry with
initial space $\pi(E_{s^*})H$. Doing the same to $\nu_s\nu_s^*$,
we conclude that $\pi(E_s)H$ is the final space of $\nu_s$.\\
Let us split the demonstration that $\nu$ is a homomorphism in two
cases. Firstly take $h\in\pi(E_{t^*s^*})H$. Then
$h=\rho(a\delta_e)(k)$, $a\in E_{t^*s^*}$ and $k\in H$. Let
$\{u_\lambda\}$ an approximately identity of $E_s$. Using the
first part of the demonstration of that $\nu$ is well defined we
have:
\begin{equation*}
\begin{split}
\nu_s\nu_t(h)&=\nu_s\nu_t(\rho(a\delta_e)(k))=\nu_s\rho(\beta_t(a)\delta_t)(k)=
\D\lim_\lambda\rho(u_\lambda\delta_s)\rho(\beta_t(a)\delta_t)(k)=\\
&=\D\lim_\lambda\rho(\beta_s(\beta_{s^*}(u_\lambda)\beta_t(a))\delta_{st})(k)=
\D\lim_\lambda\rho(u_\lambda\beta_s(\beta_t(a))\delta_{st})(k)=\\
&=\rho(\beta_{st}(a)\delta_{st})(k)=\nu_{st}(h).
\end{split}
\end{equation*}
Now let $h\in(\pi(E_{t^*s^*})H)^\perp$. We have that
$\nu_{st}(h)=0$. Let us show that $\nu_s\nu_t(h)=0$. Take
$\{u_\lambda\}$ an approximately identity of $E_s$ and
$\{u_\gamma\}$ of $E_t$. Well,
$\beta_s(\beta_{s^*}(u_\lambda)u_\gamma)\in\beta_s(E_{s^*}\cap
E_t)=E_{st}$ and by the Cohen - Hewitt Factorization Theorem,
$\beta_s(\beta_{s^*}(u_\lambda)u_\gamma)=xy, x,y\in E_{st}$. By
hypothesis:
\begin{equation*}
\langle\rho(\beta_{t^*s^*}(y)\delta_e)(h),H\rangle=\langle
h,\rho(\beta_{t^*s^*}(y)\delta_e)H\rangle=\langle
h,\pi(\beta_{t^*s^*}(y))H\rangle=0,
\end{equation*}
that implies $\rho(\beta_{t^*s^*}(y)\delta_e)(h)=0$. As
\begin{equation*}
\rho(\beta_s(\beta_{s^*}(u_\lambda)u_\gamma)\delta_{st})(h)=\rho(xy\delta_{st})(h)=
\rho(x\delta_{st})\rho(\beta_{t^*s^*}(y)\delta_e)(h)=0,
\end{equation*}
taking $\{u_\lambda\}\subset E_s$ and $\{u_\omega\}\subset
E_{s^*}$ their approximately identities, follows that:
\begin{equation*}
\nu_s\nu_t(h)=\D\lim_\lambda\rho(u_\lambda\delta_s)\D\lim_\omega\rho(u_\omega\delta_t)(h)
=\D\lim_{\lambda,\omega}\rho(\beta_s(\beta_{s^*}(u_\lambda)u_\omega)\delta_{st})(h)=0.
\end{equation*}
So $\nu_{st}=\nu_s\nu_t$ and $\nu$ is a
homomorphism.\\
For last, we need to prove the covariance condition, that is, that
$\nu_s\pi(a)\nu_{s^*}=\pi(\beta_s(a))$. Let $a\in E_{s^*}$ and
$\{u_\lambda\},\{u_\gamma\}$ approximately identities of $E_s$.
Then:
\begin{equation*}
\begin{split}
\nu_s\pi(a)\nu_{s^*}&=
\D\lim_\lambda\rho(u_\lambda\delta_s)\rho(a\delta_e)\D\lim_\gamma\rho(\beta_{s^*}(u_\gamma)\delta_{s^*})=
\D\lim_{\lambda,\gamma}\rho(u_\lambda\delta_s)\rho(a\beta_{s^*}(u_\gamma)\delta_{s^*})=\\
&=\D\lim_{\lambda,\gamma}\rho(\beta_s(\beta_{s^*}(u_\lambda)a\beta_{s^*}(u_\gamma))\delta_{ss^*})=
\D\lim_{\lambda,\gamma}\rho(u_\lambda\beta_s(a)u_\gamma\delta_{ss^*})=\\
&=\rho(\beta_s(a)\delta_{ss^*})=\rho(\beta_s(a)\delta_e)=\pi(\beta_s(a)),
\end{split}
\end{equation*}
because $\rho|_N\equiv0$. \vspace{-0.5cm}
 \begin{flushright}

  $\blacksquare$

  \end{flushright}

By the previous Theorem we have, for $x\in L$:
\begin{equation*}
\begin{split}
\|x\|_c&=\sup\{\|(\pi\times\nu)(x)\|:
(\pi,\nu,H)\in\hbox{CovRep}(A,S,\beta)\}=\\
&=\sup\{\|\rho(x)\|:\rho\hbox{ representation of }L\hbox{ equals zero in }N\}=\\
&=\sup\{\|\rho(x)\|:\rho\hbox{ representation of
}A\times_\beta^aS\}=\rho_s(x).
\end{split}
\end{equation*}

Then:

\begin{theo}
The definition of Partial Crossed Product by an Action of Inverse
Semigroup in a C$^*$-Algebra that we present is equivalent to that
introduced by Sieben in \cite{sieben}.
\end{theo}

 \begin{flushright}

  $\blacksquare$

  \end{flushright}

\begin{center}
References
\end{center}
\vspace{-1.7cm}
\renewcommand{\refname}{}


\begin{thebibliography}{99}
\bibitem{dokuchaexel} M. Dokuchaev and R. Exel, \emph{Associativity
of Crossed Products by Partial Actions, Enveloping Actions and
Partial Representations}. Trans. American Math. Soc. \textbf{357}
(2005), 1931-1952.

\bibitem{exel} R. Exel, \emph{Partial actions of groups and actions
of inverse semigroups}. Proc. Am. Math. Soc. \textbf{126} (1998),
3481-3494.
%
%\bibitem{mccl} K. McClanahan, \emph{K-theory for partial crossed
%products by discrete groups}. Journal of Funct. Anal. \textbf{130}
%(1995), 77-117.

\bibitem{hewitt} E. Hewitt and K. A. Ross, \emph{Abstract Harmonic Analysis
II}, Springer-Verlag, 1970.

\bibitem{lawson} M. V. Lawson, \emph{Inverse Semigroups - The Theory of Partial
Symmetries}, World Scientific Publishing Co. Pte. Ltd., 1998.

\bibitem{sieben} N. Sieben, \emph{C$^*$ crossed products by partial
actions and actions of inverse semigroups}. Journal Austral. Math.
Soc. Ser. A \textbf{63} (1997), 32-46.

\bibitem{sunder} V. S. Sunder, \emph{Functional Analysis - Spectral Theory},
Birkhäuser Verlag, 1998.
\end{thebibliography}
\end{document}